\pdfoutput=1
\documentclass[12pt]{article}
\usepackage{amsmath,amsfonts,amssymb,color,a4,epsfig,graphics}
\usepackage[latin1]{inputenc}

\newcommand{\R}{{\mathbb{R}}}
\newcommand{\C}{{\mathbb{C}}}
\newcommand{\Z}{{\mathbb{Z}}}
\newcommand{\N}{{\mathbb{N}}}

\def\ha{\frac{1}{2}}

\def\pa{\partial}
\def\ra{\rightarrow}
\def\preuve{\begin{proof}}
\def\ga{\alpha}

\def\gd{\delta}
\def\ge{\varepsilon}

\def\gg{\gamma}

\def\gl{\lambda}
\def\go{\omega}

\def\gr{\rho}
\def\gs{\sigma}

\def\es{essentially self-adjoint}
\newtheorem{defi}{Definition}[section]

\newtheorem{lemm}{Lemma}[section]
\newtheorem{prop}{Proposition}[section]
\newtheorem{rem}{Remark}[section]
\newtheorem{coro}{Corollary}[section]
\newtheorem{theo}{Theorem}[section]

\newtheorem{ques}{Question}[section]
\newenvironment{demo}{\noindent {\it Proof.--}
      \begin{quotation}\noindent}{\end{quotation}\hfill$\square $}

\renewcommand{\footnote}[1]{\footnotetext{#1}}

\begin{document}

\footnote{ \textbf{Math Subject Classification (2000):} 05C63, 05C50, 05C12, 35J10, 47B25.}
\title{Essential self-adjointness for combinatorial Schr\"odinger
operators II- Metrically non complete graphs}
\author{Yves Colin de Verdi\`ere\footnote{Grenoble University,
Institut Fourier,
 Unit{\'e} mixte
 de recherche CNRS-UJF 5582,
 BP 74, 38402-Saint Martin d'H\`eres Cedex (France);
{\tt yves.colin-de-verdiere@ujf-grenoble.fr};
{\tt http://www-fourier.ujf-grenoble.fr/$\sim $ycolver/}
}\\ Nabila Torki-Hamza
\footnote{Unit\'e de recherches Math\'ematiques et Applications, 
05/UR/15-02, Facult\'e des Sciences de Bizerte de l'Universit\'e de Carthage, 
7021-Bizerte; ISIG-K, Universit\'e de Kairouan, 3100-Kairouan; Tunisie;
{\tt nabila.torki-hamza@fsb.rnu.tn};
{\tt natorki@yahoo.fr}}
\\ Fran{\c c}oise Truc\footnote{Grenoble University, Institut Fourier,
Unit{\'e} mixte
 de recherche CNRS-UJF 5582,
 BP 74, 38402-Saint Martin d'H\`eres Cedex (France);
{\tt francoise.truc@ujf-grenoble.fr};
{\tt http://www-fourier.ujf-grenoble.fr/$\sim$trucfr/}
}}
\date{March 3, 2013}
\maketitle


\begin{abstract}

We consider weighted graphs, we equip them with a metric structure given
by a weighted distance, and we discuss essential self-adjointness for weighted graph Laplacians
and Schr\"odinger operators in the metrically non complete case.
\end{abstract}

\section{Introduction}
This paper is a continuation of \cite{To} which contains some statements about
 essential self-adjointness of Schr\"odinger operators on graphs.
In \cite{To}, it was proved that on any metrically complete weighted graph
with bounded degree, the Laplacian is essentially self-adjoint
and the same holds for the Schr\"{o}dinger operator provided the
associated quadratic form is bounded from below. These results remind
those in the context of Riemannian manifold in \cite{Ol}
and also in \cite{B-M-S}, \cite{Shu1}, \cite{Shu2}.
There are many recent independent researches in locally finite graphs investigating essential
self-adjointness (see \cite{Jor}, \cite{Go-Sch}, \cite{Ma}), and relations between stochastic completeness
and essential self-adjointness ( see \cite{We},  \cite{Woj2} as well as the thesis \cite{Woj1}).
Similar results have been extended for arbitrary regular Dirichlet forms on discrete sets
in \cite{Ke-Le-2} which is mostly a survey of the original
article \cite{Ke-Le-1}. More recently the paper \cite{Hu}
is devoted to the stability of stochastic incompleteness, in almost the same setup as in \cite{Ke-Le-1}.
\newline
Here, we will investigate essential self-adjointness mainly on
metrically non complete locally finite graphs.
\newline
Let us recall that a weighted graph $G$
is a generalization of an electrical network where the set
of vertices and the set of edges are respectively weighted with positive
functions $\go$ and $c$. For any given positive function $p$,
a weighted distance $d_p$ can be defined on $G$. Thus we have the usual notion of
completeness for $G$ as a metric space.
\newline
The main result of Section \ref{non-esa}
states that the weighted graph Laplacian $\Delta_{\go,c}~$
(see the definition (\ref{equ:deflap}) below) is not essentially
self-adjoint if the graph is of finite volume and metrically non complete
(here the metric $d_p~$ is defined using the weights
$p_{x,y}=c_{x,y}^{-\ha}~$).
The proof is derived from the existence of the  solution
for a Dirichlet problem at infinity, established in Section \ref{dirichlet}.
\newline
In Section \ref{noncomp},  we establish some conditions implying
essential self-adjointness. More precisely, defining
the metric $d_p$ with respect to the weights $p_{x,y} ~$ given by
$p_{x,y} =(\min \{ \go_x, \go_y\})c_{x,y}^{-\ha}~,$
 and addressing the case of metrically non complete graphs,
 we get the essential self-adjointness of $\Delta_{\go,c}+W$
under the assumptions that the potential $W$ is bounded from below by $ N/2D^2~ $, where $N$ is the maximal
degree and $D$ the distance to the boundary, and that the graph has a regularity property.
We use for this result a technical tool
deduced from Agmon-type estimates and
inspired by the nice paper \cite{Nen}, see also \cite{Col-Tr}.
\newline
We discuss in Section \ref{treelike} the case of star-like graphs.
Under some assumptions on $a$, we prove that for any potential $W$,
$\Delta _{1,a}+ W$ is \es ~  using an extension
of Weyl's theory to the discrete case. In the particular case
of the graph $\N$, the same result had been proved
in \cite{Ber} (p.504) in the context of Jacobi matrices.
We give some examples in Subsection \ref{examples} to illustrate
the links between the previous results. Moreover we establish the sharpness of the conditions
of Theorem \ref{essod}.
\newline
The last Section is devoted to  Appendix A dealing with Weyl's limit point-limit circle
criteria (see \cite{RS}) in the discrete case as well as in the continuous case,
and to Appendix B including the unitary equivalence between Laplacians and Schr\"{o}dinger
operators \cite{To} used repeatedly in Subsection \ref{examples}.
\newline
Let us start with some  definitions.
\newline
$G =(V,E)$ will denote an {\it infinite  graph}, with $V=V(G)$
the set of {\it vertices} and $E=E(G)$ the set of {\it edges}. We write
 $ x\sim
y$ for  $\{x,y\}\in E$.
\newline
The graph $G$ is always assumed to be \emph{locally
finite}, that is  any $x\in V$ has a finite number
 of neighbors, which we call  the \emph{degree}
(or valency) of $x$.  If the degree is bounded
independently of $x$ in $V$, we say that the graph $G$ is of \emph{bounded
 degree}.
\newline
 The space of real functions on the graph $G$ is denoted
$$C(V)=\left\lbrace f:V\longrightarrow \R\right\rbrace $$
and $C_0(V)$ is the subspace of functions with finite support.\\
We consider, for any weight $\go:V\longrightarrow ]0,+\infty[$, the
space
  $$l_\go^2(V)=\lbrace f\in C(V);\; \underset{x\in
V}{\sum}\go_x^2 f^2(x) <\infty \rbrace.$$
It is a Hilbert
space when equipped with  the inner product:
\[ \langle f,g \rangle_{ l_{\omega}^2 } =\sum_{x\in V}
 \omega_x^2 f\left( x\right).g\left( x\right)~.\]
For any  $\go:V\longrightarrow ]0,+\infty[$, and $c:E\longrightarrow ]0,+\infty[$,
the \emph{weighted graph Laplacian} $\Delta_{\go,c}$ on the graph
$G$ weighted by the \emph{conductance} $c$ on the edges and by the \emph{weigth} $\go$
on the vertices, is defined by:
 \begin{equation} \label{equ:deflap}
 \left( \Delta_{\go,c}f\right)\left( x\right)
 =\dfrac{1}{\go_{x}^{2}} \sum_{y\sim x}
 c_{ x,y}\left( f\left( x\right) -f\left(y \right) \right)
 \end{equation}
for any $f\in C(V)$ and any $x\in V .$
 If $\go \equiv 1$, we have
$$\left( \Delta_{1,a}f\right)\left( x\right) = \sum_{y\sim x}
 a_{ x,y}\left( f\left( x\right)-f\left(y \right) \right). $$
\begin{defi} \label{def:dist}
Let $p:E\longrightarrow ]0, +\infty [$  be given,
the  {\rm weighted distance} $d_p (\leq +\infty)$
 on the weighted graph $G$ is defined by
$$d_p(x,y)=\inf_{\gamma \in \Gamma_{x,y}}
 L(\gamma )$$
where $\Gamma_{x,y} $ is the set
 of the paths $\gamma : x_1=x, x_2, \cdots, x_{n}=y $
 from $x$ to $y$. The length $L(\gamma)$ is computed as
the sum of the $p$-weights for the edges of the path $\gamma $:
$$L(\gamma)=\sum_{1\leq i\leq n} p_{x_i,x_{i+1}}~.$$
In particular, if $x$ and $y$ are in distinct connected components
of $G$, $d_p(x,y)=\infty$.
We say that the metric space $(G,d_p)$ is complete when every Cauchy sequence
 of vertices has a limit in $V$.
\end{defi}

\begin{defi}\label{def:compl}
We denote by ${\hat V}$ the metric {\rm completion}
 of $(G, d_p)$ and by  $ V_{\infty} = \widehat{V}
\setminus V $ the metric {\rm boundary} of $V$.
\end{defi}
\begin{defi} \label{def:ends}
If $G$ is a non finite graph and $G_0
$ a finite sub-graph of $G$, the  {\bf ends of $G$ relatively to
$G_0 $} are
the  non finite connected components of $G\setminus G_0 $.
\end{defi}
\section{The Dirichlet problem at infinity}\label{dirichlet}
We will use in this section the  distance $d_p$  defined
 using the weights
$p_{x,y}=c_{x,y}^{-\ha}$.
Let us consider the quadratic form
\[ Q(f)= \sum_{\{x,y\} \in E} c_{x,y}(f(x)-f(y))^2 +
 \sum _{x\in V}\go_x^2 f(x)^2 ~,\]
which is formally associated to the operator $\Delta _{\go, c}+{\rm
  Id}$ on $l^2_\go$.
We will need the following result
 which is close  to lemma 2.5 in \cite{Jo-Pe-2}:
\begin{lemm}\label{lemm:lip}
For any $f:V \ra \R $ so that $Q(f)<\infty$
and for any $a,b  \in V$,
we have
$$|f(a)-f(b)| \leq \sqrt{Q(f)}d_p(a,b).$$
\end{lemm}
\begin{demo}
For any $\{x,y\} \in E$,
$|f(x)-f(y)| \leq \sqrt{Q(f)}/\sqrt{c_{x,y}}$.
For any path $\gg$ from $a$ to $b$,  defined by the
vertices  $x_1=a, x_2, \cdots, x_{n}=b $,
we have
$|f(a)-f(b)|\leq  \sqrt{Q(f)}L(\gg ) $.
Taking the infimum of the righthandside with respect to $\gg$
we get the result.
\end{demo}
\begin{rem}
Lemma \ref{lemm:lip} implies that any function $f$ with
$Q(f)<\infty $ extends to $\hat{V}$ as a Lipschitz function $\hat{f}$.
We will denote by $f_\infty $ the restriction
of  $\hat{f}$ to $V_\infty. $
\end{rem}

\begin{theo}\label{functionF} Let us assume that $(V,d_p)$ is non complete.
Let $f:V\ra \R $ with $Q(f) <\infty $, then there exists a
continuous
function $F:\hat{V}\ra \R $ which satisfies both conditions:
\begin{description}
  \item (i)  $(F - {f})_{\infty}=0 $
  \item (ii) $(\Delta _{\go,c} +1)(F_{|V})=0 $.
\end{description}
Moreover, such an $F$ satisfies $Q(F)<\infty $ and $F\in l^2_\go $.

If $\hat{V}$ is compact, such an $F$ is unique.
\end{theo}
\begin{demo}
We will denote by $A_f$ the affine space of continuous functions $G:
\hat{V}\ra \R $ which satisfy $Q(G)<\infty $
and $({G}-{f})_{\infty}=0$.

$Q $ is lower semi-continuous for the pointwise convergence
on $V$ as defined by  $Q=\sup Q_\ga  $ with
$Q_\ga  (f)=$ sum of a finite number of terms in $Q$.

Let $Q_0:= \inf_{G\in A_f}Q(G)$ and
 $G_n$ be  a corresponding minimizing sequence.
The  $G_n$'s are equicontinuous and pointwise bounded.
From Ascoli's Theorem,
 this implies the existence of a locally uniformly  convergent
subsequence
$G_{n_k} \ra F$.
Using semi-continuity, we have $Q(F)=Q_0$.

If $x\in V$ and $\gd_x $ is the Dirac function at the
vertex $x$, we have
\[ \frac{d}{dt}_{\mid t=0} Q(F+ t\gd _x)=
2\go_{x}^2[( \Delta _{\go,c}+1 )F  (x)]~\]
and this is equal to $0$, because $F$ is a minimum
of $Q$ restricted to $A_f$.

Uniqueness is proved using a maximum principle: let us assume that
there
exists a non zero continuous  $F$ with $F_\infty =0 $, then,
changing, if necessary, $F$ into $-F$,
there exists $x_0 \in V$ with $ F(x_0)=\max_{x\in V} F(x)>0 $.
The identity (ii) evaluated at the vertex $x_0$ gives a contradiction.
\end{demo}

\section{Not essentially self-adjoint
Laplacians }\label{non-esa}

\begin{theo}\label{theo:nonESA}
 Let $\Delta _{\go, c}$ be a weighted graph Laplacian and assume
the following conditions:

\begin{description}
  \item (i) $(G,d_p)$ with $p_{x,y}=c_{x,y}^{-\ha}$ is NON complete,
  \item (ii) there exists a function $f:V\ra \R $ with
$Q(f) <\infty $ and
$f_{\infty }\ne  0 .$
\end{description}
Then  $\Delta _{\go, c}$ is not essentially self-adjoint.
\end{theo}
\begin{demo}
Because $\Delta _{\go,c}$ is $\geq 0 $ on $C_0(V)$,
 it is enough (see Theorem X.26 \cite{RS}) to
build a non zero function $F:V \ra \R $
which is in $l^2_\go (V) $ and satisfies
 \begin{equation} \label{equ:delta+1}
(\Delta _{\go,c}+ 1)F=0 ~.
\end{equation}
 The function $F$ given by Theorem \ref{functionF} will be the
 solution of equation (\ref{equ:delta+1})
 the limit of which at infinity is $f_{\infty }$.
\end{demo}

\begin{rem}\label{nonco} The assumptions of Theorem \ref{theo:nonESA}
are satisfied if $(G, d_p)$ is non complete and
$\sum \go_y^2 <\infty $: it is enough to take $f\equiv 1$.

They are  already satisfied if $G $
 has a non complete   ``end'' of finite volume.
\end{rem}
\begin{rem}
Theorem \ref{theo:nonESA} is not valid for the Riemannian Laplacian:
if $X$ is a closed Riemannian manifold of dimension $\geq 4 $,
$x_0 \in X$  and
$Y=X\setminus {x_0 }$, the Laplace operator
on $Y$ is  essentially self-adjoint (see \cite{Col1})
 and $Y $ has finite volume.
\end{rem}
\begin{ques}
 In Theorem \ref{theo:nonESA}, what is the deficiency index
of $\Delta _{\go,c}$ in terms of the geometry
of the weighted graph?\end{ques}


\section{Schr\"odinger operators for
metrically non complete graphs}\label{noncomp}

We now discuss essential selfadjointness for  Schr\"{o}dinger operators of the type $H=\Delta _{\go, c} +W $ on a graph $G$ in
the following setup: we define $\alpha_{x,y}= \min \{ \go_x, \go_y\}$
 and we assume  that $(G,d_p)$, with  $p_{x,y} = \dfrac{\alpha_{x,y}}{\sqrt{c_{x,y}}}~,$
 is {\rm non complete} as a metric space.
It means that there exist Cauchy sequences of vertices
without limit in the set $V$.
We will assume that $G$ is of bounded degree, and we denote
 the upper bound by $N$. We will need also to assume a regularity property for $(G,d_p)$.
\subsection{Regularity property for metric graphs}
\begin{defi}\label{dist}
 For a vertex $x\in V$,
 we denote by $D(x)$ the distance to the boundary  $ V_{\infty}$ defined by
\begin{equation}\label{dist} D(x)= \inf_{ z \in  V_{\infty}} d_p(x,z).\end{equation}
\end{defi}

\begin{lemm}\label{D}
We have, for any edge $\{ x,y \}$,
\begin{equation}|D(x)- D(y)| \leq d_{p}(x,y) \leq \dfrac{\min \{ \go_x, \go_y\}}{\sqrt{c_{x,y}}}
\ .
\end{equation}
\end{lemm}
\begin{defi}If $A $ is a subset of $ V$, the boundary  $\pa
A$ of $A$ is  
 the set of the vertices 
$x\in  A $ so that there exists $y \in V \setminus A $
 with $\{x,y\} \in E$.\end{defi}
\begin{defi}\label{ikx} Let $\ge >0$ be given and   $ X_{\ge  }$ be defined by
    \[ X_{\ge  }=\{ x\in V ~|~ D(x) \geq \ge \}~.
\] 
We say that the graph $(G,d_p)$ is {\rm regular} if, for any sufficiently small $\ge$, any bounded
subset 
of $\pa X_\ge $ (for the metric $d_p$)  is finite.
\end{defi}

The main property  of the regular graphs that we will use is:
\begin{prop}\label{finsup} If $(G,d_p)$ is regular, then closed and bounded subsets of
$X_{\varepsilon  }$ are finite.
\end{prop}
\begin{demo} Let $A$ be a closed and bounded subset of $X_{\varepsilon}$.

Let $\Gamma=(V', E')$ be the graph with  vertices $V' =A$
and edges $E'= E'_1 \cup E'_2$ defined as follows:

\medskip

\noindent (i) edge $\{ x,y \} $ is in $E'_1$ iff   $\{ x,y \}\in E$ and  $x,\, y\in A$,

\medskip

\noindent (ii) edge $\{ x,y \} $ is in $E'_2$ iff  $\{x,y\}\notin E$ and $x,\,y\in A\cap \partial X_{\varepsilon}$.

\medskip

Note that the edges $E'_1$ already exist in $G=(V,E)$, while the set $E'_2$ represents new edges.

The graph $\Gamma$ is equipped with the following weight:
\[
 p'_{x,y}= \left\{
 \begin{array}{cc}
            p_{x,y}, & \{x,y\} \in E'_1,\\
            d_p(x,y), & \{x,y\} \in E'_2.
            \end{array}
            \right.
 \]
Clearly, the distances $d_{p'} $ and $d_p$  coincide on $A $.  If $(x_n)_{n\in \N }$ is a Cauchy sequence
for $(A, d_p)$,  it has  a limit in $\hat{V}$ which
lies in $A $ because $A$ is closed in $\hat{V}$.
Hence, the graph $(\Gamma,d_{p'})$ is complete. Since $A\cap \partial X_{\ge}$ is a bounded subset of $\partial X_{\varepsilon}$, by regularity property it 
follows that $A\cap \partial X_{\varepsilon}$ is finite; hence, the set $E'_2$ is finite. This, together with local finiteness of $G=(V,E)$,
 shows that $\Gamma=(V',E')$ is locally finite. Since $A $ is closed and bounded with respect to the metric $d_p$, it is also closed and 
bounded with respect to $d_{p'}$. Thus, we  can apply the discrete version of Hopf-Rinow Theorem for the graph $(\Gamma,d_{p'})$ (see \cite{HKMW}, Theorem
  A1) to conclude that $A$ is finite.
\end{demo} 

\subsubsection{Trees are regular}
Let $T$ be a locally finite tree with the weight $p$.
 Let us  choose a root $x_0$
of
$T$.
\begin{defi} A {\rm ray} is a maximal simple path in $T$ starting from $x_0$.
\end{defi}
Let us start with two Lemmas:
\begin{lemm} The Cauchy boundary of the weighted tree $(T,d_p)$
  identifies
with the set of rays of finite length with an infinite number
of vertices. More precisely, if  the ray $(x_0,x_1,\cdots,x_n,\cdots )$
has a finite length, the sequence $(x_n) $ is a Cauchy sequence.
The map $\Phi $ which  associates the corresponding
point in $T_\infty $ to this Cauchy sequence is a bijection. 
\end{lemm}
\begin{demo} The map $\Phi $ is injective, because if $(x_n)$
and $(y_l)$ are two rays of finite length the distance $d_p(x_n, y_l)$
is bounded from below  by a strictly positive number as soon
as $n$ or $l$ is large enough.

Similarly, if $(x_n)$ is a Cauchy sequence, then $x_n$ has to be on a
single  ray  for
$n$ large enough. Hence $\Phi $ is surjective. 
\end{demo}
 
\begin{lemm}\label{finiteray}
If all rays of a  tree $T$ have a finite number of vertices, $V(T)$ is
finite.
\end{lemm}
\begin{demo}
Let us prove the Lemma by contradiction. Let us denote by $|x|\in \N$
 the
combinatorial distance from $x_0$ to $x$.
If $V(T)$ is not finite, we can build by induction,
 using the local finiteness of $T$, a ray   $(x_n)_{n\in \N}$ with 
$|x_n|=n $, so that the number of vertices of the sub-tree rooted at
$x_n$ and contained in $\{ |x|\geq n \}$ is infinite.
This ray is infinite, hence a contradiction.
\end{demo}

\begin{prop} If $T$ is a tree, then for any choice of the weight $p$,
$(T,d_p)$ is regular.\end{prop}

\begin{demo}
 Let us consider a ball $B(x_0,R)$ 
in $(T,d_p)$. The rays of infinite length are going out of 
$B(x_0,R)$ after a finite number of steps. 
The set $X_\ge $ is obtained by removing on each infinite  ray of finite length
all vertices with $|x|$ large enough.
 Hence, the set $X_\ge \cap
B(x_0,R)$
is the set of vertices of a tree all of whose rays have a finite
number
of vertices, and we can apply Lemma \ref{finiteray}.
 \end{demo}

\begin{coro}
If the first Betti number $b_1(G)$ is finite, then $(G,d_p)$ is regular.
\end{coro}
\begin{demo}
 Removing a finite number of edges, we get a tree  and the
  sets 
$X_\ge$ are the same for both graphs
if $\ge $ is small enough, as well as the bounded sets of $V$.
  \end{demo}

\subsubsection{An example with a finite Cauchy boundary not satisfying
the regularity condition}

Let us consider the following graph $G=(V,E)$:
$V=\N \times \Z $,
\[ E=\{ \{ (n,k), (n+1,k)\}| n\in \N, k\in \Z \}
\cup \{ \{ (n,k), (n,k+1)\}| n\in \N, k\in \Z \} ~,\]
and choose the weight $p$ defined as follows:
 \[ p_{ \{ (n,k), (n+1,k)\}}= 2^{-(n+1)} ,~   p_{ \{ (n,k), (n,k+1)\}}
 =1/n ~.\]

The distance $d_p$ satisfies, for $(n,k)\ne (n',k')$:
\[ 2^{-(1 +\min (n,n'))} \leq
   d_p ( (n,k),(n',k')) \leq 2^{-n}+ 2^{-n'}\leq 1 ~.\]
The lower bound is just the min of the weights of the edges from
one of the ends of the path; the upper bound is the 
limit of the lengths of the curves $\gamma _N$ as $N \ra \infty $
where $\gg_N$ is described as  follows, assuming  that 
$k'\geq k$ and $N> \max(n,n')$: $\gg_N$  starts  from $(n,k)$, follows
 $(n+1,k),(n+2,k),\cdots ,(N,k)$, then 
   $(N,k),(N,k+1),\cdots , (N,k')$
and finally
$(N,k'),(N-1,k'),\cdots, (n',k')$.

Then $(G,d_p)$ is bounded;  a sequence $x_l=(n_l,k_l)$
 is a Cauchy sequence if and only if $n_l
\ra \infty $; more precisely the diameter of the set 
$X_N:=\{ (n,k)~|~n\geq N \}$ is less than $2^{-(N-1)}$.  
 The Cauchy boundary of $(G,d_p)$ is reduced to a single point.
 The function $D$ is given by
\[ D( (n,k))=2^{-n}~.\]
Hence the graph $(G,d_p)$ is not regular.

\subsection{Agmon-type estimates}

\begin{lemm}\label{nen}

Let $v,~f \in C_{0}(V)$ be real valued
and assume  $Hv=0$.
Then
 \begin{equation}\label{ute}
\langle   fv\ , H (fv) \rangle_{l^2_{\go}}   \ =\
  \sum_{ \{ x,y \}  \in E} c_{x,y} v(x)
  v(y) (f(x)-f(y))^2~.
\end{equation}

\end{lemm}

\begin{demo}
In the case of positive $v$ this type of formula is known as ground state
transform (see \cite{Hea-Kel} and references within). A particular case of this computation (for operators of the
type $\Delta _{1, a} +W $) can be found in \cite{To}, let us recall the proof for the reader's convenience:
$$\langle   fv\ , H (fv) \rangle_{l^2_{\go}}
= \ \sum_{x\in V}f(x)v(x)\left(
 \sum_{y\sim x}\ c_{x,y}(f(x)-f(y))v(y)\right)$$
where we used the fact that  $Hv(x)=0$.
An edge $\{x,y\}$ contributes to the  sum twice.
 The total contribution is
$$f(x)v(x)\ c_{x,y}(f(x)-f(y))v(y) +f(y)v(y)
 c_{y,x}(f(y)-f(x)) v(x)$$ so
$$\langle fv\ , H (fv) \rangle_{l^2_{\go}}=\sum_{\{x, y\} \in E}
 c_{x,y}(f(x)-f(y))\left( f(x)v(x) v(y) -f(y)v(y)
 v(x) \right)\ . $$
\end{demo}

\begin{theo}\label{Horo}
Assume  that $(G,d_p)$, with  $p_{xy} = \dfrac{\alpha_{x,y}}{\sqrt{c_{x,y}}}~,$
 is a {\rm non complete regular} graph. Let $v$ be a  solution of $(H-\gl)v=0$.
 Assume that $v$ belongs to $l^2_{\go}(V )$
and that there exists
  a constant $c>0$ such that, for all $u \in C_0(V)$,
\begin{equation}\label{boun}
\langle u | (H-\gl) u \rangle_{l^2_{\go}}  \geq
 \dfrac{N}{2}\sum_{x\in V}\max \left(\dfrac{1}{D(x)^2},1
\right) \go_x^2  |u(x)|^2  +  c \|u\|_{l^2_{\go}}^2~,
\end{equation}
then $v\equiv 0$.
\end{theo}

\begin{demo}
This theorem is based on Lemma \ref{nen}
 applied to $H-\gl$.
Let us consider  $\rho $ satisfying
$0< \rho < \dfrac{1}{2}$.
For any $\ge >0$, we define the  function $f_\ge:V \rightarrow \R $
 by $f_\ge =F_\ge(D) $ where  $D$ denotes the distance associated to the  metric $d_p$ , as in  (\ref{dist}), and $F_\ge: \R^+ \rightarrow \R $ is the
continuous piecewise affine
  function defined by
\[ F_\ge(u)= \left\{
\begin{array}{l}
0  {\rm ~ for~}  u\leq \ge   \\
\gr  (u-\ge)/(\gr  - \ge  )   {\rm ~ for~} \ge  \leq u \leq  \gr  \\
u    {\rm ~ for~} \gr  \leq u \leq  1  \\
1   {\rm ~ for~} 1 \leq u 
\end{array}
\right.
\]

Let us fix a vertex $x_0$. For any $\ga >0$, we define also the  function $g_\ga:V \rightarrow \R $
 by $g_\ga =G_\ga(d_p(x_0,.)) $ where  $G_\ga: \R^+ \rightarrow \R $ is the
continuous piecewise affine
  function defined by
\[ G_\ga(u)= \left\{
\begin{array}{l}
1  {\rm ~ for~}  u\leq 1/\ga   \\
-\ga  u + 2  {\rm ~ for~} 1/\ga  \leq u \leq  2/\ga  \\
0   {\rm ~ for~ }  u \geq 2/\ga
\end{array}
\right.
\]

Let  $E_{\ge , \ga} $  be the set of the vertices defined by 

\begin{equation}\label{support}
 E_{\ge , \ga} = \{x \in V ~|~ \ge \leq D(x) \  {\rm and}  \ d_p(x_0,x) \leq 2/\ga \}\ .
\end{equation}
Note that the support of $f_{\varepsilon}g_{\alpha}$ is contained in $E_{\varepsilon, \alpha}$. Additionally, note that $E_{\varepsilon, \alpha}$ is a closed 
and bounded subset of $X_{\varepsilon}$, where $X_{\varepsilon}$ is as in Definition \ref{ikx}. Hence, $E_{\varepsilon, \alpha}$ 
is finite by Proposition \ref{finsup}. Therefore, the function $f_{\varepsilon}g_{\alpha}$ is finitely supported.
Observe moreover that
$$|f_{\varepsilon}(x) g_{\alpha}(x)-f_{\varepsilon}(y)g_{\alpha}(y)|\leq |f_{\varepsilon}(x)| |g_{\alpha}(x)-g_{\alpha}(y)| +|g_{\alpha}(y)|
|f_{\varepsilon}(x)-f_{\varepsilon}(y)|$$
$$\leq \dfrac{\gr}{\gr -\ge}|D(x) -D(y)|   +  \ga |d_p(x_0,x)-d_p(x_0,y)|
$$
so using Lemma \ref{D}  we get that $f_\ge g_\ga$ is $(\dfrac{\gr}{\gr -\ge} + \ga)-$Lipshitz with respect to the metric $d_p$.

We can apply Lemma \ref{nen} to
the finite-supported function $f_\ge g_\ga$, and using  the inequalities
$$  v(x)v(y) \leq \ha ( v(x)^2+v(y)^2 )~,$$
we  get that the right hand side  of (\ref{ute}) is bounded as follows

$$
\langle   f_\ge  g_\ga v | (H-\gl) (f_\ge  g_\ga v) \rangle_{l_{\go}^2}  \leq
\ha \left(\dfrac{\gr}{\gr -\ge }+\ga \right)^2\sum _{x\in V } v(x)^2  \Phi_{\ge , \ga}  (x) ~,$$
with
 $$ \Phi_{\ge , \ga}  (x) =\sum _{y\sim x } c_{x,y} d_p(x,y)^2
   \leq N \go_x^2  ~ $$

where the first inequality uses the fact that
$f_\ge g_\ga$ is $(\dfrac{\gr}{\gr -\ge} + \ga)-$Lipshitz with respect to the metric $d_p$,
 and the second inequality is a direct consequence of the choice of the weights $p_{x,y}$ (see Lemma \ref{D}).
This implies

\begin{equation} \label{equ:l2hs}
\langle   f_\ge g_\ga v | (H-\gl) (f_\ge g_\ga v) \rangle_{l_{\go}^2}  \leq \frac{N}{2}
\left(\dfrac{\gr}{\gr -\ge }+\ga \right)^2 \| v \| _{l^2_\go}^2 ~.
\end{equation}

On the other hand, due to assumption (\ref{boun}) the left hand side
of  (\ref{ute}) is bounded from below as follows:

\begin{equation} \label{equ:r2hs}
\langle   f_\ge g_\ga v | (H-\gl) (f_\ge g_\ga v) \rangle_{l_{\go}^2}
\geq \frac{N}{2} \sum _{E_{\ge , \ga}} \go_x^2  v(x)^2 +c \| f_\ge g_\ga v
\|_{l_{\go}^2}^2~,
\end{equation}
where $E_{\ge , \ga} $ is as in (\ref{support}).

Putting together  (\ref{equ:l2hs}) and (\ref{equ:r2hs}) we get

\begin{equation} \label{aest}
 \frac{N}{2} \sum _{E_{\ge , \ga}}\go_x^2  v(x)^2 +c \| f_\ge g_\ga v
\|_{l_{\go}^2}^2  \leq \frac{N}{2} \left(\dfrac{\gr}{\gr -\ge }+\ga \right)^2  \| v \| _{l^2_\go}^2 \ .
\end{equation}

Then we do $\ga \ra 0 $. After that, we do also $\ge \ra 0 $ . The last step is to take the limit  $\gr \ra 0 $,
and then we get that  $v\equiv 0 $.

\end{demo}

\begin{rem}

The previous result is inspired by a nice idea from \cite{Nen}, so following  the terminology
of \cite{Nen}  we call
 Agmon-type estimates Equation (\ref{aest}) .
\end{rem}
\subsection{Essential self-adjointness}

\begin{theo}\label{essod}
 Consider  the Schr\"{o}dinger operator $H=\Delta _{\go, c} +W $ on a graph $G$, define $\alpha_{x,y}= \min \{ \go_x, \go_y\}$
 and assume  that $(G,d_p)$, with  $p_{xy} = \dfrac{\alpha_{x,y}}{\sqrt{c_{x,y}}}~,$
 is a {\rm non complete regular} graph.
For a vertex $x\in V$,
 we denote by $D(x)$ the distance from $x$ to the boundary  $ V_{\infty}$.
 We assume the following conditions:
 \begin{description}
   \item (i) $G$ is of bounded degree and we denote the upper  bound by $N$,
   \item (ii) there exists $M < \infty $ so that
\begin{equation}\label{buno}
 \forall x \in V, ~ W(x) \geq \ \dfrac{N}{2D(x)^2}-M~.
 \end{equation}
 \end{description}
Then the Schr\"{o}dinger operator $ H$ is essentially self-adjoint.
\end{theo}

\begin{rem} In the particular case when $\sum \go_x^2< \infty $, the Laplacian
$H=\Delta _{\go, c} $ does not satisfy the assumption (\ref{buno}) so this result
is coherent with Theorem \ref{theo:nonESA}.
\end{rem}

\begin{rem}
The exponent of $D(x)$ in (\ref{buno}) is sharp. In fact, one can find a potential $W$ such that
$W(x) \geq \ \dfrac{k}{D(x)^2}~$ where $k<\dfrac{N}{2}$ and weights
$\go$ and $c$ such that $H=\Delta_{\go, c}+W$ is non essentially self-adjoint. See Example \ref{exam2}~.
\end{rem}

\begin{rem}\label{com} In the case where  $\go \equiv 1$
the result is an immediate consequence of  \cite{Ke-Le-1} (Theorem 5).
\end{rem}
\begin{demo}
We have, for any $u \in C_0(V)$

$$\langle u | H u \rangle_{l_{\go}^2} \geq \sum_{x\in V} W(x)\go_x^2 |u(x)|^2, $$
so using assumption (\ref{buno}) we get:
$$\langle u | (H-\gl) u \rangle_{l^2_{\go}}  -\frac{N}{2}\sum_{x\in V} \frac{1}{D(x)^2}
\go_x^2   |u(x)|^2
\geq   \sum_{x\in V}\
 -(M+\gl) \|u\|_{l^2_{\go}}^2 ~ . $$

Then choosing for example $$ \gl
= -M-1$$ we get the
 inequality (\ref{boun})
with $c=1$, and the proof follows from Theorem \ref{Horo}.
\end{demo}

\section{Schr\"odinger operators on
 ``star-like'' graphs}\label{treelike}

\subsection{Introduction}

\begin{defi}
The graph $\N$ is the  graph
 defined by $V=\{ 0,1,2,\cdots \} $
and $E=\{ \{n,n+1\}~|~n=0,1,\cdots \}$.
\end{defi}

\begin{defi}
We will call an infinite graph $G=(V,E)$  {\rm star-like}
if there exists a finite sub-graph $G_0$ of $G$
so that $G\setminus G_0$ is the union of a finite number
of disjoint copies $G_\ga $  of the graph $\N$ (the {\rm ends} of $G$
relatively to $G_0$ according to Definition \ref{def:ends}).
\end{defi}

For example, the graph $\Z$, defined similarly to
$\N$, is star-like.

Let us consider a Laplace operator $L=\Delta _{1,a}$ on $G$.
On each end  $G_\ga$ of $G$, $L$ will be given by
\[ L^\ga  f_n \ = \ - a_{n,n+1}^\ga f_{n+1} +
 ( a_{n-1,n}^\ga + a_{n,n+1}^\ga ) f_{n}- a_{n-1,n}^\ga f_{n-1}\
,\]
where the $a_{n-1,n}^\ga$'s are $>0$.
If $W:V\ra \R $, we will consider Schr\"{o}dinger
 operators $H$ on $C_0(G)$
defined by
 $H = \Delta_{1,a}  + W\ .$

\begin{lemm}\label{lemm:ends}
 Let $G_0 $ be a finite sub-graph of $G$.
The operator $H=\Delta_{1,a}  + W$ on $G$
is essentially self-adjoint if and only if it is
 essentially self-adjoint on each end of $G$ relatively to $G_0$.
 More precisely, the deficiency indices $n_\pm$
 are the sum
of the corresponding deficiency indices of the ends.
\end{lemm}

We will need the following Lemma which is a consequence of Kato-Rellich
Theorem, see \cite{Go-Sch}, Proposition 2.1:

\begin{lemm}\label{lemm:index}
If $A $ and $B$ are 2 symmetric
operators with the same domains and $R=B-A$ is bounded,
then the deficiency indices of $A$ and $B$ are the same.
\end{lemm}
\begin{demo}
We give here an alternative  proof to this result.
Let us define, for $t\in \R$,  $A_t=A+tR $ so that $A_0=A$
and $A_1=B$.
The domains of the   closures of the $A_t$'s  co\"incide:
the ``graph-norms''
$ \| A_tu \|_{l^2}  + \| u \| _{l^2 }$
 are equivalent.
The domains of the adjoints
co\"incide too.
Let $K= D(A^\star)/ D( \bar{A})$ and
$Q_t (u,v)=-i\left(
\langle A_t^\star u |v \rangle -\langle u|  (A_t^\star v \rangle
\right)$
which is well defined on $K$. We know that these bounded Hermitian
forms are non degenerate on $K$ with the
graph norm and continuous w.r. to $t$.
Hence the Morse index $n_- (t)$
is locally constant: take
a decomposition $K=K_+ \oplus K_- $ where $q=Q_{t_0} $
satisfies $q_{|K_+}\geq C >0 $ and $q_{|K_-}\leq -C <0 $.
\end{demo}

Using Lemma \ref{lemm:index}, we can prove
Lemma \ref{lemm:ends}:

\begin{demo}
We will consider the operator $H_{\rm red}$
where we replace the entries $a_{x,y}$ of $H$
with $\{ x, y \} \in E(G_0 )$ by $0$.
The claim of the Lemma is clear for $H_{\rm red}$ because
it is the direct orthogonal sum of the Schr\"{o}dinger operators of the ends
and a finite rank $l^2-$bounded matrix.
We can then use Lemma \ref{lemm:index} because
$H-H_{\rm red}$ is bounded.
\end{demo}

\begin{rem}\label{starg}
It follows from Lemma \ref{lemm:ends} that, concerning
 essential self-adjointness questions  for star-like graphs,  it is enough to
work on the graph $\N$.
We have
$$(Hf)_0=- a_{0,1}f_{1} +  a_{0,1} f_{0}+ W_0 f_0 .$$
This implies that the space of solutions of $(H-\gl)u=0$
on $\N$ is of dimension $1$ and any solution so that $f_0$ vanishes
is $\equiv 0$.
We will consider also solutions ``near infinity'',
 i.e. $(f_n)_{ n\geq 0}$ satisfies
 $((H-\gl)f)_n=0$ for $n\geq 1$; this space is
of dimension $2$.
\end{rem}

\subsection{Main result}

 It is known (\cite{Dod})
that $H=\Delta_{1,a}+W$ is essentially self-adjoint provided $\Delta_{1,a}$ is
bounded as an operator on $l^2(G)$ and $W$ bounded from below.
For star-like graphs, we have the following result, which holds for {\it any} potential $W~:$

\begin{theo}\label{ess}
 If $G$ is star-like and if for each end $G_\alpha$,
  \begin{equation}\label{prop}
 1/a_{n-1,n}^\ga \notin l^1(\N)
\end{equation} then
 $H=\Delta_{1,a}  + W$ with domain $C_0(V)$ is essentially self-adjoint for any potential $W$.
\end{theo}

\begin{rem}
The  condition (\ref{prop}) is sufficient but not necessary. See Example \ref{exam2}~.
\end{rem}

\begin{demo}
Due to Remark \ref{starg} we only have to prove
the following

\begin{theo}\label{essn}
If
  \begin{equation}\label{propo}
\frac{1}{a_{n-1,n}} \notin l^1(\N) ~,
\end{equation}
the Schr\"{o}dinger operator ${H}= \Delta_{1,a}  + W$ with
 domain $C_0(\N)$ is essentially self-adjoint for any potential $W$~.
\end{theo}
This result is contained in the book \cite{Ber} (p. 504).
 We propose here a short proof, obtained  by contradiction
 using Corollary \ref{Weyl} which is an analog
of Weyl's limit point-limit circle criteria in the discrete case.

Let us consider an operator $\Delta_{1,a}$ such
that (\ref{prop}) is fulfilled.
We assume that any sequence  $u$, such that $(H-i)u =0$
near infinity, is  in $l^2(\N)$.
In particular, there exists a basis $f,g$ of solutions
 of ${ (H-i)}f =0$  with $f \in l^2(\N)$ and $g \in l^2(\N)$.

We have
$$-a_{n,n+1}\ f_{n+1} + ( a_{n-1,n}+ a_{n,n+1} +(W_n-i))
 \  f_{n}- a_{n-1,n}\ f_{n-1}\ =\ 0\ ,$$
and the same holds for $g$.
The Wronskian of $f$ and $g$ is the sequence
 ${\cal W}_n=f_n g_{n-1}- f_{n-1}g_n $.
We have, for any $n \in \N$:
$${\cal W}_{n+1} =
 \frac{a_{n-1,n}}{ a_{n,n+1}}\ {\cal W} _n\ ,$$
which implies
 \[ {\cal W}_n\ =\ \frac{a_{0,1}}{ a_{n-1,n}}{\cal W}_1 \quad
 .\]
But since the Wronskian is in $l^1(\N)$ according to the assumption
 that  $f$ and $g$ are in $l^2(\N)$, we get a contradiction with
the hypothesis (\ref{prop}).

\end{demo}


\subsection{ Examples of   Schr\"odinger
operators }\label{examples}

\subsubsection{Example 1}

 Let us consider the Laplacian $\Delta _{\go, c}$ on $\N$,
 with, $ \forall n >0$,  $c_{n-1,n}= n^3$ and,  $ \forall n \geq 0 $,
$\omega_n=\dfrac{1}{ n+1}$. Since $\sum \go_n^2 < \infty$
and $\sum c_{n-1,n}^{-1/2} < \infty$ we deduce from Theorem
\ref{theo:nonESA} (due to Rem \ref{nonco} ) that $\Delta _{\go, c}$
is not essentially self-adjoint.

Applying a result of \cite{To} (see Proposition \ref{equiv} in Appendix B)
we get that this Laplacian is unitarily equivalent to
the Schr\"odinger operator  $H= \Delta_{1,a}  + W$
 with $a_{n-1,n}= \dfrac{c_{n-1,n}}{\go_{n-1}\go_n}\sim n^5$ and
 $$W_n = \dfrac{1}{\go_n}
 \left[c_{n,n+1}\left(\dfrac{1}{\go_n}-\dfrac{1}{\go_{n+1}}\right) +
c_{n-1,n}\left(\dfrac{1}{\go_n}-\dfrac{1}{\go_{n-1}}\right) \right]\sim -3n^3 ~,$$
which is therefore not essentially self-adjoint.

 According to Theorem \ref{ess}, such an operator must verify
 $ \dfrac{1}{a_{n-1,n}} \in l^1(\N)$, which is indeed the case.
\subsubsection{Example 2:  Discretization of a Schr\"odinger
operator on  $\R^+$}\label{exam2}

 Let us consider the Schr\"odinger operator on  $]0,+\infty [$ defined
 on smooth compactly supported functions by
$Lf:= -f" + \dfrac{A}{x^2} f$.
This operator is essentially self-adjoint if and only if $A > 3/4$ (see \cite{RS} theorem X 10).
We discretize this operator in the following way:
let us consider the graph $\Gamma = (V,E)$ resulting of
 the following dyadic subdivision of the interval $(0,1) $:
the vertices are the $x_n =2^{-n}$ and the edges
 are the pairs $\{2^{-n}, 2^{-n+1}\}$
which correspond to the intervals
 $[2^{-n}, 2^{-n+1} ]$  of length $\go_n^2 =2^{-n}$.

Then we define,  for any
 \[ f \in l_\go^2(V)=\lbrace f\in C(V)~|~
\sum_{n\in \N} 2^{-n} f_n^2 <
 +\infty \rbrace\]  where we set $f= (f_n)$,
 the quadratic form
  $$Q(f)= \sum_{n\in \N}  2^{-n}\left[\left( \dfrac{f_{n+1}-f_n}{2^{-n}}\right)^2
 + A 2^{2n} f_n^2 \right] \ .$$
According to the previous definitions and if we set $c_{n,n+1}= 2^{n}$,
this quadratic form is associated to the Schr\"odinger operator
$H=\Delta _{\go, c} + W$ on $\N$ with the potential $W_n:=A 2^{2n}$.

 Let us set $a_{n,n+1}=\dfrac{c_{n,n+1}}{\go_n\go_{n+1}}=2^{2n+\ha }$.
Applying  Proposition \ref{equiv} we get that $H$ is unitarily equivalent
to
$$\widehat{H}= \Delta_{1,a}  + \widehat{W} + W$$
with
 $$ \widehat{W}_n = \dfrac{1}{\go_n}
 \left[c_{n,n+1}\left(\dfrac{1}{\go_n}-\dfrac{1}{\go_{n+1}}\right) +
c_{n-1,n}\left(\dfrac{1}{\go_n}-\dfrac{1}{\go_{n-1}}\right) \right] =
  2^{2n}\left(\dfrac{3}{2}-\dfrac{5\sqrt{2}}{4}\right)\ .$$
We have $\widehat{H}=\Delta_{1,a}  + (A-A_0) 4^n$ with
   $A_0 =\dfrac{5\sqrt{2}}{4} -\dfrac{3}{2}~(>0)$.
The metric graph  $(\N, d_p)$ with $p_{n,n+1}=
a_{n,n+1}^{-1/2}$ is non complete.
The solutions $u$ of  $Hu= 0$ verify
$$4 u_{n+1} - \left( 5 +2\sqrt{2}(A -A_0) \right) u_n +  u_{n-1} =0\ .$$
The solutions are generated by $\ga_1  ^n $
and $\ga_2 ^n $
where $\ga_1$ and $\ga _2$ are the roots of
\[ 4\ga^2 -\left( 5 +2\sqrt{2} (A -A_0) \right)\ga +1=0 ~.\]
We have
$|\ga_1|<1$ and $|\ga _2|<1 $ if and only if $A_0-\dfrac{5}{\sqrt{2}}< A< A_0 $.

Using Proposition \ref{prop:dyn}, with $d=2$
and $U_n=\left( \begin{array}{c} u_n \\ u_{n-1} \end{array}  \right)
$, we get, for any $\gl \in \C$, the exponential decay of all solutions near infinity
of  $(\widehat{H}-\gl)u=0$ if  $A_0-\dfrac{5}{\sqrt{2}}< A< A_0 $, and the existence of
a solution of $(\widehat{H}-\gl)u=0$ with exponential growth in the case when $ A > \dfrac{5\sqrt{2}}{4} -\dfrac{3}{2} $
or   $A < -\dfrac{5 \sqrt{2}}{4} -\dfrac{3}{2}$.

 Hence  (by Corollary \ref{Weyl}) we get the following result:

\begin{prop}\label{pr}

\begin{enumerate}

\item If  $-\dfrac{5\sqrt{2}}{4} -\dfrac{3}{2} < A< \dfrac{5\sqrt{2}}{4} -\dfrac{3}{2} $,
then the discretized operator $H$ is not essentially self-adjoint.

\item If $ A > \dfrac{5\sqrt{2}}{4} -\dfrac{3}{2} ~(\star)~$\quad
or   $A < -\dfrac{5 \sqrt{2}}{4} -\dfrac{3}{2}$, then $H$ is essentially self-adjoint.
\end{enumerate}\end{prop}

From this result we can deduce several informations:
\begin{enumerate}
  \item The condition $ (\star)$ is analogous to the condition
$A>3/4$  in the continuous case.
  \item Proposition \ref{pr} implies that for $A =0$ the operator $H = \Delta_{\go,c} $
is not essentially self-adjoint, which is a result predicted by Theorem  \ref{theo:nonESA}.
  \item This gives examples of essentially self-adjoint operators with $1/a_n \in l^{1} $.
  \item Sharpness of the assumption (\ref{buno}) in  Theorem \ref{essod}

In this context, the distance $d_p$ is associated to
$$p_{x,y}= \dfrac{\alpha_{x,y}}{\sqrt{c_{x,y}}}$$
with $\alpha_{x,y}= \min \{ \go_x, \go_y\}$ so
we get $$D(n) = \sum _{p\geq n}
\dfrac{\alpha_{p,p+1}}{\sqrt{c_{p,p+1}
}}=  \sum _{p\geq n}
\left(\dfrac{2^{-p-1}}{2^{p}}\right)^{1/2}\ = \ 2^{-\ha} 2^{-n} 2$$ so $$\dfrac{1}{D(n)^2 } = 2^{2n-1}\ .$$

If the operator $H=\Delta_{\go,c}  + A 4^n$ satisfies the
assumption (\ref{buno}),
 then \\ $A>\dfrac{1}{2}$ which involves condition $(\star)$,
 since $\dfrac{1}{2}>\dfrac{5\sqrt{2}}{4} -\dfrac{3}{2}$, so Theorem
 \ref{essod} is  coherent with
proposition \ref{pr}.
Moreover the operator $H=\Delta_{\go,c}  + A 4^n$ with $A =\dfrac{5\sqrt{2}}{4} -\dfrac{3}{2} $ is not essentially self-adjoint,
which implies that the estimate (\ref{buno}) on the growth of the potential in  Theorem \ref{essod} is sharp.
\end{enumerate}


\subsubsection{Example 3}

 Let us consider the Laplacian $\Delta _{\go, c}$ on $\N$,
 where the coefficients verify
 $c_{n-1,n}= n^\gamma $
with $\gamma >2$  and
$\omega_n=(n+1)^{-\beta }$ with $\beta >\dfrac{1}{2} $.
Since $\sum \go_n^2 < \infty$ and $\sum c_{n-1,n}^{-1/2} < \infty$ we deduce from Theorem
\ref{theo:nonESA} (due to Remark \ref{nonco} )
 that $\Delta _{\go, c}$ is not essentially self-adjoint.

Applying one more time Proposition \ref{equiv}, we see that this operator is unitarily equivalent to
the Schr\"odinger operator $H= \Delta_{1,a}  + W$, with
 $a_{n-1,n} \sim n^{\gamma +2 \beta } $ and the potential
 $W_n \sim -\beta(\beta +\gamma -1)n^{2\beta +\gamma -2} $,
 which is therefore also not essentially self-adjoint.
 We emphasize that $W$ is not bounded from below, which is predicted
 in \cite{To}, Theorem 3.2.

 Furthermore, according to Theorem \ref{ess},
 such an operator must verify the condition $ \dfrac{1}{a_{n-1,n} }\in l^{1}(\N)$,
 which is indeed the case.
Following the terminology of the previous sections,
 it means  the non completeness of $(\N,d_p)$ with the weights $p_{n-1,n}=
a_{n-1,n}^{-1/2}~.$

\subsubsection{Example 4}

Let us consider the Laplacian $H= \Delta _{\go, c}$
 on a spherically homogeneous rooted tree $G=(V,E)~$
 (see \cite{Bre} and references within).
For any vertex $x$, we  denote by ${\delta}(x)$ the
 distance from $x$ to the root $0$
and define
$ \go_x=2^{-\delta(x)}$, and
 $c_{x,y}= 2^{\delta(x)}$, for any $y\sim x$ so that ${\delta}(y)=n+1~$.
We assume that the graph $G$ has a uniform degree $N+1$.

 Let us set $a_{x,y}=\dfrac{c_{x,y}}{\go_x\go_y}$.
 We have $a_{x,y}=2^{3n+1 } $ for any edge ${x,y}$,
 so that ${\delta}(x)=n$ and $ {\delta}(y)=n+1$.
 Then, due to Proposition \ref{equiv}, the operator $H$ is unitarily equivalent
to
$$\widehat{H}= \Delta_{1,a} + W$$
with
 $$W(x) = 2^{3n}\left(-N +\dfrac{1}{4}\right) $$
 for any $ x $ such that $ {\delta}(x) =n \ .$

The radial solutions $u$ of $Hu= 0$ can be seen as sequences
$(u_n)$ which satisfy the equation:
$$ -2N u_{n+1} + \left( N +\dfrac{1}{2}\right) u_n -\dfrac{1}{4} u_{n-1} =0\ .$$

The solutions are generated by $\ga_1  ^n $
and $\ga_2 ^n $
where $\ga_1$ and $\ga _2$ are the roots of
\[ \ga^2 -\left( \dfrac{1}{2} -\dfrac{1}{4N} \right)\ga +\dfrac{1}{8N}=0 ~.\]
We have
$|\ga_1|<1$ and $|\ga _2|<1 $ for any $N>0$.

The radial solutions of $(\widehat{H}-\gl)u=0$
satisfy
$$ -2N u_{n+1} + \left( N +\dfrac{1}{2}\right) u_n -\dfrac{1}{4}
u_{n-1} =2N \gl 2^{-(3n+1) } u_n
\ .$$
Using Proposition \ref{prop:dyn}, with $d=2$
and $U_n=\left( \begin{array}{c} u_n \\ u_{n-1} \end{array}  \right)
$, we get the exponential decay of all solutions near infinity
of  $(\widehat{H}-\gl)u=0$.

 Hence  (by Corollary \ref{Weyl}) we get the following result:
\begin{prop}\label{pro}
For any  $ N \geq 1$ $H$ is not essentially self-adjoint .

\end{prop}

\begin{rem}

We have $$\sum_{x}\go_x^2 =\sum_{n}\go_n^2  N^n = \sum_{n}(\dfrac{N}{4})^n \ .$$
If $N< 4$, then
$\sum_{x}\go_x^2 < \infty$ so Theorem  \ref{theo:nonESA} can also be applied to get
the result since the graph is non complete with respect
to the metric $d_p$ , with $p_{x,y}= c_{x,y}^{-\ha}$.

\end{rem}


\section{Appendix A: Weyl's  ``limit point-limit circle'' criteria}
\label{appendix}
\subsection{The discrete case}\label{appendiscrete}
The goal of this section is to prove the discrete version
of the Weyl's ``limit point-limit circle'' criterium.
Our presentation is simpler than the classical
presentation for the continuous case (see \cite{RS}, Appendix to
section
X.1).

Let us consider the Hilbert space
${\cal H}:=l^2(\N, \C^N )$
and the formally symmetric differential
operator
$P$ defined by
\[ Pf(0)= P_{0,0} f(0)+ P_{0,1}f(1),~\forall l\geq 1,
  Pf(l)=P_{l,l-1}f(l-1)+ P_{l,l} f(l)+ P_{l,l+1} f(l+1 ) ~\]
where
\begin{enumerate}
\item $\forall l\geq 1,~P_{l-1,l}^\star = P_{l,l-1}$
\item  $\forall l\geq 0,~P_{l,l}^\star = P_{l,l}$
\item $ \forall l\geq 0,~P_{l,l+1} $ is invertible
\item $  \exists M\in  \R$ so that
for any $f \in C_0 (\N, \C^N)$,
$Q_P(f)=\langle Pf ~|~f \rangle \geq -M \| f\| ^2 $.
\end{enumerate}

Let us define the subspace ${\cal E} $
of $ {\cal H}$ as the set of $l^2$ sequences $f$ so that,
for all $l\geq 1$,
$(P-i)f(l)=0 $; the space
   ${\cal E} $  is isomorphic to the space of germs at
infinity
of $l^2$ solutions of $(P-i)f=0$.
 Assumption 3. implies that  $\dim {\cal E}\leq 2N$.
Let us denote by ${\cal K}= \ker (P-i)\cap l^2$ and
 consider the following sequence
\begin{equation} \label{sequ} 0 \ra {\cal K}
 \ra  {\cal E} \ra \C^N \ra 0 ~,\end{equation}
where the non trivial arrow is given by
$f \ra (P-i)f(0)$.
We have the
\begin{theo}
The sequence (\ref{sequ}) is exact
and the deficiency indices
$n_\pm = \dim {\cal K}$ of $P$  are given by
$n_\pm = \dim  {\cal E}-N $.
\end{theo}
\begin{demo}
 Assumption 4. implies (using Corollary of Theorem X.1
 in \cite{RS}) that the deficiency indices are equal.
The only non trivial point is to prove that the arrow
$p: {\cal E}\ra \C^N $ is surjective.
Let us consider $\tilde{P}$ a self-adjoint extension of $P$
which exists because $n_+=n_-$.
Let us consider the map
$\rho : \C^N \ra  {\cal E}$ defined
by
$\rho (x)= (\tilde{P}-i)^{-1} (x,0,0,\cdots ) $.
Then $p\circ \rho = {\rm Id}_{\C^N}$.
\end{demo}

\begin{coro}\label{Weyl}
The Schr\"{o}dinger operator $H= \Delta_{1,a}  + W$
 defined on $C_0(\N)$ is essentially self-adjoint  if and only if
 there exists a sequence $u$ such that $(H-i)u =0$
near infinity  (i.e. $((H-i)u)_n=0$ for $n$ large enough)
which is not in $l^2(\N)$.
\end{coro}
\subsection{Asymptotic behavior of perturbed
hyperbolic iterations}

In order to apply Corollary \ref{Weyl}, the following results will
be useful
\begin{prop}\label{prop:dyn}
Let us consider the following linear dynamical system on $\C^d$:
\begin{equation} \label{equ:dyn}
\forall n\geq 0,~ U_{n+1}=A U_n + R(n) U_n \end{equation}
where
\begin{enumerate}
\item $A$ is hyperbolic:  all eigenvalues $\gl_j $ of $A$
 satisfy $|\gl_j |\ne  1 $
\item  $R(n) \ra 0 $ as $n \ra \infty $.
\end{enumerate}
Then
\begin{itemize}
\item {\bf case A:} If all eigenvalues  $\gl_j $ of $A$
 satisfy $|\gl_j |< 1 $,  all solutions $(U_n)$
 of Equation (\ref{equ:dyn}) are
exponentially
decaying.

\item {\bf case B:} If  $m$  eigenvalues satisfy  $|\gl_j |>  1 $,
then there exists an  $m$-dimensional vector space $F$
of solutions  of Equation (\ref{equ:dyn}) whose non-zero vectors have
exponential growth.
\end{itemize}
\end{prop}

\begin{demo}

{\bf Case A:}
 There exists a norm $\| . \| $ on $\C^d$ so that
the operator norm of $A$ satisfies
$\| A \| =k < 1 $.
For $n$ large enough, we have $\| A + R(n) \| \leq k' < 1$.
The conclusion follows.

{\bf Case B:} There exists a splitting $\C^d=Y\oplus Z$, denoted
$x=y+z $, with $\dim Y =m$,
stable by $A$, norms on $Y$ and $Z$ and
2  constants $\mu < 1 <\sigma $,  so that
\[ \forall y\in Y, \| Ay \| \geq \sigma \| y \| ~,\]
\[ \forall z\in Z, \| Az\| \leq \mu \| z \| ~.\]
Let us choose $\ge > 0 $ so that
$1 < \gs  -2\ge $ and $N$ so that
$\| R(n) \| \leq \ge $ for $n\geq N$.
We have, for $n\geq N $,
\[ \| y_{n+1}\|  \geq \gs \| y_n \| - \ge (\| y_n \| +\| z_n \| ),
  \| z_{n+1}\|  \leq \mu  \| z_n \| + \ge (\| y_n \| +\| z_n \| )~,\]
so that
\[  \| y_{n+1}\|- \| z_{n+1}\|\geq (\gs -2\ge) ( \| y_{n}\|- \|
z_{n}\|)~.\]
Any solution which satisfies $\| y_N \| > \| z_N \| $
will have exponential growth. Take for $F$  the space of
solutions for which $z_N=0$.

 \end{demo}

\subsection{The continuous case}

A similar method works for the continuous case.
Let
$H=-\dfrac{d^2}{dx^2 } + A(x) $ be a system of differential operators
where $A(x) $
is Hermitian for every $x$ and is continuous on
$[0,a[ $ as a function of $x$. The differential operator
$H$ is $L^2$-symmetric on the Dirichlet domain
 \[D= C_0^\infty( [0,a[,\C^N) \cap \{ u ~|~ u(0)=0 \}~.\]
We denote $H_D$ the closure of $(H,D)$.
Let us assume that $n_+(H_D)=n_-(H_D)$ which is true for example
if $A $ is bounded from below or if $A$ is real-valued.
Then
\begin{theo}
If ${\cal E}$ is the space of solutions $ u$
of the differential equation \\ $(H-i)u=0$ which are  $L^2 $ near $a$,
then
$n_\pm (H_D)= \dim {\cal E} -N $.
\end{theo}
\begin{demo}
Let us consider the sequence
\begin{equation} \label{sequ-cont} 0 \ra \ker (H_D -i)
 \ra  {\cal E} \ra \C^N \ra 0 ~,\end{equation}
where the only non trivial arrow is given
by $u \ra u(0)$.
This sequence is exact: we have only to prove the surjectivity
of the non trivial arrow.
Let $\tilde{H}$ be a self-adjoint extension
of $H_D$ and  $\chi \in C_0^\infty ([0,a[,\R) $ with $\chi (0)=1$.
For any  $X\in \C^N$, let us consider
\[ u=\chi X - (\tilde{H}-i)^{-1}\left( (H-i) (\chi X )\right)~.\]
Then $(H-i)u=0 $, $u(0)=V$ and $u$ is $L^2$ near $a$.
\end{demo}

\section{Appendix B: Unitary equivalence between Laplacians and Schr\"{o}dinger operators }
\label{appendixB}

In this section, we recall the following results
(see \cite{To} Proposition 2.1 and Theorem 5.1):
the first one states that a Laplacian is always
unitarily equivalent to a Schr\"{o}dinger operator,
and the second result asserts that a Schr\"{o}dinger
operator with a strictly positive quadratic form is unitarily
equivalent to a Laplacian.

For a weighted graph $G$ by the weight $\go$ on its vertices, let
\begin{equation*}
U_\omega :~l_{\omega}^2 \left( V\right)\longrightarrow l^2 \left( V\right)
\end{equation*}
the unitary operator defined by
\begin{equation*}
U_\omega\left( f\right) =\omega f~.
\end{equation*}
This operator preserves the set of functions on $V$ with finite support.

\begin{prop}\label{equiv}
 The operator
\begin{equation*}
\widehat{\Delta}= U_\omega ~\Delta_{\omega,c}~U_\omega ^{-1}~,
\end{equation*}
is a Schr\"odinger operator on $G$~. More precisely:
\begin{equation*}
\widehat{\Delta}=\Delta_{1,a}+W
\end{equation*}
where $a$ is a strictly positive weight on $E$ given by:
\begin{equation*}
a_{ x,y}=\dfrac{c_{ x,y}}{\omega_x\omega_y}
\end{equation*}
and the potential $W:V\longrightarrow {\mathbb{R}}$ is given by:
\begin{equation*}
W=-\dfrac{1}{\omega}~\Delta_{1,a}\omega~.
\end{equation*}
\end{prop}

The following Theorem uses the existence of a strictly positive
harmonic function (see \cite{To},  section 4).

\begin{theo}
Let $P$ a Schr\"odinger operator on a graph $G$~.
We assume that $\langle Pf,f\rangle_{l^{2}} > 0$ for any function $%
f$ in $C_0\left( V\right)\setminus\lbrace0\rbrace $~. Then there exist
 weights: $\go$ on $V$ and  $c$
on $E$ such that $P$ is unitarily equivalent to the
Laplacian $\Delta_{\omega,c}$ on the graph $G$~.
\end{theo}

\section*{Acknowledgments}

We would like to thank Ognjen Milatovic, who pointed out the existence of some gap
in the proof of  Theorem \ref{Horo} of our early version of the paper. In order to
remove this gap, we have made an additional assumption on the graph, which is a regularity
property , and which allows us
 to prove the finiteness of the support of the cut off functions in Theorem \ref{Horo}.

The second author is greatly indebted to the research unity
"Math\'ematiques et Applications" (05/UR/15-02) of
Facult\'e des Sciences de Bizerte (Tunisie) for the financial support,
and would like to present special thanks to Institut Fourier where this work was carried on.
\newline
Thanks to D. Lenz for giving notes on some references.
\newline
All the authors would like to thank the reviewers for their comments and 
especially Reviewer 1 for
careful reading, numerous remarks, useful suggestions and valuable references.


\end{document}